\documentclass[11pt]{article}

\usepackage{graphicx} 
\usepackage{amsmath,amssymb,amsthm,mathtools}
\usepackage{xcolor}
\usepackage{hyperref}
\usepackage[capitalise]{cleveref}
\usepackage{enumerate}
\usepackage{caption}
\usepackage{subcaption}
\usepackage{tikz}
\usepackage{placeins}
\usetikzlibrary{calc}
\usetikzlibrary{arrows.meta}

\usepackage[letterpaper,twoside,outer=1.3in,vmargin=1.3in,]{geometry}

\newcommand{\ignore}[1]{}

\title{Strong orientation of a connected graph for a crossing family
}
\author{Ahmad Abdi \and Mahsa Dalirrooyfard \and Meike Neuwohner}

\date{\today}

\theoremstyle{theorem}
\newtheorem{theorem}{Theorem}

\newtheorem{conjecture}[theorem]{Conjecture}

\begin{document}

\maketitle

\begin{abstract}
Given a connected graph $G=(V,E)$ and a crossing family $\mathcal{C}$ over ground set $V$ such that $|\delta_G(U)|\geq 2$ for every $U\in \mathcal{C}$, we prove there exists a strong orientation of $G$ for $\mathcal{C}$, i.e., an orientation of $G$ such that each set in $\mathcal{C}$ has at least one outgoing and at least one incoming arc. This implies the main conjecture in Chudnovsky et al.\ (\emph{Disjoint dijoins}. Journal of Combinatorial Theory, Series B, 120:18--35, 2016). In particular, in every minimal counterexample to the Edmonds-Giles conjecture where the minimum weight of a dicut is $2$, the arcs of nonzero weight must be disconnected.
\end{abstract}

\section{Introduction}

A set collection $\mathcal{C}\subseteq 2^V$ is called a \emph{crossing family} over the ground set $V$ if $U\cap W,U\cup W\in \mathcal{C}$ for all sets $U,W\in\mathcal{C}$ that \emph{cross}, i.e., $U\cap W\neq \emptyset$ and $U\cup W\neq V$. Given a graph $G=(V,E)$ and a crossing family $\mathcal{C}$ over ground set $V$, a \emph{strong orientation of $G$ for $\mathcal{C}$} is an orientation $\overrightarrow{G}$ of the edges of $G$ such that $\delta^+_{\overrightarrow{G}}(U), \delta^-_{\overrightarrow{G}}(U) \neq \emptyset$ for every $U\in \mathcal{C}$. In this note, we prove the following theorem by using tools from integer programming and combinatorial optimization. 

\begin{theorem}\label{thm:strong_orientation}
Let $G=(V,E)$ be a connected graph, and let $\mathcal{C}$ be a crossing family over ground set $V$ such that $|\delta_G(U)|\geq 2$ for all $U\in \mathcal{C}$. Then there exists a strong orientation of $G$ for $\mathcal{C}$.
\end{theorem}

\begin{figure}[t]
\begin{subfigure}[t]{0.45\textwidth}
\centering
\begin{tikzpicture}
\foreach \i in {1,...,6}
{
  \node[circle, fill, inner sep = 0pt, minimum size = 2mm] (A\i) at ({1.5*cos(\i*60)},{1.5*sin(\i*60)}) {};
  \node[circle, fill, inner sep = 0pt, minimum size = 2mm] (B\i) at ({3*cos(\i*60)},{3*sin(\i*60)}) {};
}
\foreach \i in {1,2,3,4,5,6}
{
    \draw[thick, dashed, -{Stealth}] (A\i)--(B\i);
}
\foreach \i/\j in {1/6,3/2,5/4}
{
    \draw[thick, dashed,-{Stealth}] (A\i)--(A\j);
    \draw[thick] (A\i)--(B\j);
    \draw[thick, dashed, -{Stealth}] (B\i)--(B\j);
}
\foreach \i/\j in {5/6,3/4,1/2}
{
    \draw[thick] (A\i)--(A\j);
    \draw[thick] (B\i)--(B\j);
}
\end{tikzpicture}
\subcaption{An example showing \cref{thm:strong_orientation} does not extend to disconnected graphs. $G$ is the graph on the solid edges, and $\mathcal{C}$ is the family of nonempty proper vertex subsets with no incoming dashed arc.
\label{fig:example_disconnected}}
\end{subfigure}
\hspace{1cm}
\begin{subfigure}[t]{0.45\textwidth}
\centering
\begin{tikzpicture}
\foreach \i in {1,...,6}
{
  \node[circle, fill, inner sep = 0pt, minimum size = 2mm] (A\i) at ({1.5*cos(\i*60)},{1.5*sin(\i*60)}) {};
  \node[circle, fill, inner sep = 0pt, minimum size = 2mm] (B\i) at ({3*cos(\i*60)},{3*sin(\i*60)}) {};
  \draw[thick, dashed, -{Stealth}] (A\i)--(B\i);
}
\foreach \i/\j in {1/6,3/2,5/4}
{
    \draw[thick, dashed,-{Stealth}] (A\i)--(A\j);
    \draw[thick,-{Stealth}] (A\i)--(B\j);
    \draw[thick, dashed, -{Stealth}] (B\i)--(B\j);
}
\foreach \i/\j in {5/6,3/4,1/2}
{
    \draw[thick,-{Stealth}] (A\i)--(A\j);
    \draw[thick,-{Stealth}] (B\i)--(B\j);
}
\end{tikzpicture}
    \subcaption{Schrijver's example~\cite{SCHRIJVER1980213}. Dashed arcs have a weight of $0$, solid arcs have a weight of $1$. Every dicut has weight at least $2$, but there are no disjoint dijoins among the weight-$1$ arcs.
    \label{fig:Schrijvers_example}}
\end{subfigure}
\caption{An example showing that \cref{thm:strong_orientation} does not extend to disconnected graphs (left) and Schrijver's example (right).}
\end{figure}
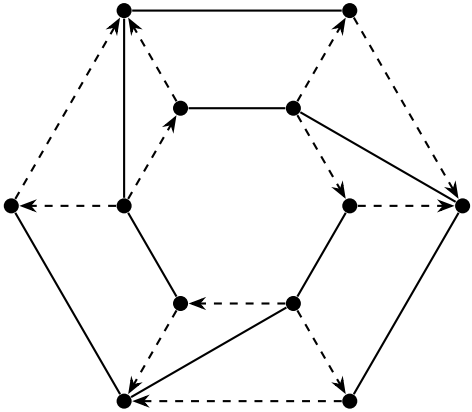
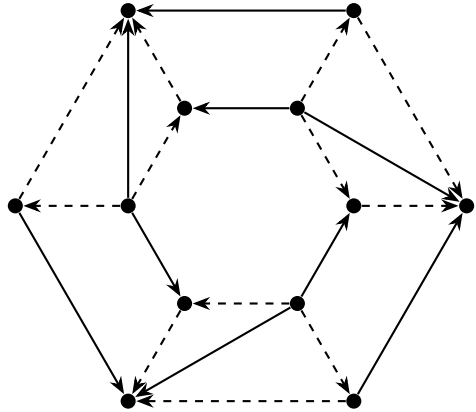

The connectivity of $G$ is crucial for this theorem, and without it the theorem would not hold. For example, let $G=(V,E)$ be the graph displayed in \cref{fig:example_disconnected} whose edges correspond to the solid lines, and let $\mathcal{C}$ be the family of all nonempty proper vertex subsets with no incoming dashed arc. (This is derived from the example in~\cite{SCHRIJVER1980213}.) Then $\mathcal{C}$ is a crossing family, $|\delta_G(U)|\geq 2$ for all $U\in \mathcal{C}$, and it can be readily checked that there is no strong orientation of $G$ for $\mathcal{C}$. 

\cref{thm:strong_orientation} proves one of the conjectures in Chudnovsky, Edwards, Kim, Scott and Seymour~\cite{CHUDNOVSKY201618} (Conjecture 2.1). In that conjecture, $G$ is a tree, and $\mathcal{C}=\{\emptyset\neq U\subsetneq V: \delta_D^-(U)=\emptyset\}$ for some digraph $D$ over vertex set $V$ (such a family is always crossing).

Our motivation for studying this problem stems from Woodall's conjecture~\cite{Woodall} and Edmonds and Giles' attempt~\cite{EDMONDS1977185} to extend this to the weighted setting. Let us elaborate.

A \emph{weighted digraph} is a pair $(D=(V,A),w)$, where $D$ is a digraph and $w\in \mathbb{Z}^A_{\ge 0}$.
Given a digraph $D=(V,A)$ and $U\subseteq V$, we define $\delta_D^+(U)\coloneqq \{(u,v)\in A\colon u\in U, v\in V\setminus U\}$ to be the set of \emph{outgoing arcs} of $U$, and we let $\delta_D^-(U)\coloneqq \{(v,u)\in A\colon u\in U, v\in V\setminus U\}$ denote the set of \emph{incoming arcs} of $U$. 
We let $\delta(U)\coloneqq \delta^+(U)\cup \delta^-(U)$.
If $D$ is clear from the context, we might omit the subscript $D$ and just write $\delta^+(U)$, $\delta^-(U)$ and $\delta(U)$, respectively.
For a vertex $v$, we write $\delta^{+/-}(v)\coloneqq \delta^{+/-}(\{v\})$ and $\delta(v)\coloneqq \delta(\{v\})$.
A \emph{dicut} in $D$ is a directed cut of the form $C\coloneqq \delta^+(U)$, where $\emptyset\neq U\subsetneq V$ such that $\delta^-(U)=\emptyset$.
A subset $J\subseteq A$ is a \emph{dijoin} if its contraction renders $G$ strongly connected. Equivalently, $J$ is a dijoin if and only if $J$ intersects every dicut.
The latter characterization implies that the maximum number $\nu$ of pairwise disjoint dijoins in $D$ is upper bounded by the minimum size $\tau$ of a dicut.
Woodall's conjecture~\cite{Woodall} asserts that $\tau=\nu$.
While it is known that if $\tau\geq 2$ then $\nu\geq 2$, Woodall's conjecture remains wide open in general, even if $\tau=3$ and the digraph is planar.

Abdi, Cornu{\'e}jols and Zlatin~\cite{Abdi2023} have shown that the arcs of $D$ can be partitioned into a dijoin and a \emph{$(\tau-1)$-dijoin}, which by definition intersects every dicut at least $\tau-1$ times; see~\cite{Abdi2024} for a short proof. Recently, by extending an idea of Shepherd and Vetta~\cite{Shepherd05}, Cornu{\'e}jols, Liu and Ravi~\cite{10.1007/978-3-031-59835-7_6} have shown that $\nu\geq \tau/6$. 

Edmonds and Giles~\cite{EDMONDS1977185} proposed a weighted version of Woodall's conjecture: Given a weighted digraph $(D=(V,A),w)$ and $\tau_w\in\mathbb{Z}_{\geq 0}$ such that the minimum weight of a dicut is at least $\tau_w$, do there exist dijoins $J_1,\dots,J_{\tau_w}$ such that for every arc $a\in A$, $|\{i\in[\tau_w]: a\in J_i\}|\leq w(a)$? 
Note that by introducing parallel arcs for arcs of nonzero weight, it suffices to consider weight vectors $w\in\{0,1\}^A$; hence, we will assume that all weights are either $0$ or $1$ in the sequel.

The Edmonds-Giles conjecture was refuted by Schrijver~\cite{SCHRIJVER1980213}, providing a counterexample for $\tau_w = 2$ (see \cref{fig:Schrijvers_example}). In this example, the weight-$1$ arcs form three distinct weakly connected components. 

The main conjecture in~\cite{CHUDNOVSKY201618} (Conjecture 1.4) is that the Edmonds-Giles conjecture is true for $\tau_w=2$ if the weight-$1$ arcs form a spanning weakly connected subdigraph. They proved this conjecture in two special cases: when $D$ is planar, and when the weight-$1$ arcs form a subdivision of a caterpillar. The proofs of these two cases, which are quite different and purely graph theoretic in nature, are technical and sensitive to extension. For example, even in the case when $D$ is planar, it was unknown prior to this work whether one could drop the `spanning' condition on the subdigraph of the weight-$1$ arcs (see step 2 in \cref{sec:proof-2}). 

In this note, we verify the conjecture of~\cite{CHUDNOVSKY201618} as a consequence of \cref{thm:strong_orientation}. A \emph{bridge} in a graph is an edge whose deletion increases the number of connected components. A graph is \emph{bridge-connected} if all its bridges, if any, belong to the same connected component. In particular, adding an isolated vertex or a disjoint $2$-edge-connected graph does not affect bridge-connectivity. A digraph is \emph{weakly bridge-connected} if its underlying undirected graph is bridge-connected. We prove the following.

\begin{theorem}\label{thm:tree_orientation}
Let $(D=(V,A),w)$ be a weighted digraph with $w\in \{0,1\}^A$, where the weight of every dicut is at least $2$. If the weight-$1$ arcs form a weakly bridge-connected subdigraph, then there exist two disjoint dijoins of $D$ contained in the weight-$1$ arcs.
\end{theorem}

This theorem extends a result of~\cite{Shepherd05} (Theorem 1.10), which proves the above in the case when every connected component of the underlying undirected graph of the weight-$1$ arcs is $2$-edge-connected.\\

The remainder of this paper is organized as follows: In \cref{sec:prelim}, we recap some basic notions and known results from combinatorial optimization and the theory of linear and integer programming that we will rely on to establish the two theorems. The proofs of \cref{thm:strong_orientation} and \cref{thm:tree_orientation} are given in \cref{sec:proof-1} and \cref{sec:proof-2}, respectively. \cref{sec:2_components} discusses the limitations of our techniques and the challenges one faces when trying to extend our approach to the case of two weakly connected components formed by the weight-$1$ arcs. We conclude with two conjectures in \cref{sec:conjecture}.

\section{Preliminaries \label{sec:prelim}}

\paragraph{Transshipments in digraphs.}

Let $D=(V,A)$ be a digraph and let $b\in \mathbb{R}^V$ with $b(V)=0$. A \emph{$b$-transshipment} is a vector $y\in \mathbb{R}^A$ with $y(\delta^+(v))-y(\delta^-(v))=b_v$ for all $v\in V$. There is a classical result providing a necessary and sufficient condition for the existence of a $b$-transshipment obeying given upper and lower bounds.

\begin{theorem}[see \cite{S2003V1}, Corollary 11.2f]\label{b-transshipment}
Let $D=(V,A)$ be a digraph, let $b\in \mathbb{R}^V$ with $b(V)=0$ and let $\ell,u\in \mathbb{R}^A$ with $\ell\leq u$. Then there exists a $b$-transshipment $y$ with $\ell\leq y\leq u$ if and only if
\begin{equation}
b(U)\leq u(\delta^+(U))-\ell(\delta^-(U)) \quad\forall U\subseteq V.   \label{eq:condition_transshipment}  
\end{equation} Moreover, if $b,\ell,u$ are integral, then $y$ can be chosen to be integral.
\end{theorem}

By setting $\ell\coloneqq-\infty$ and $u\coloneqq +\infty$, it follows from \cref{b-transshipment} that there exists a $b$-transshipment if and only if $b(U)=0$ for every $U\subseteq V$ with $\delta(U)=\emptyset$. In particular, if $D$ is weakly connected, there always exists a $b$-transshipment.

\paragraph{Crossing-submodular functions.}
Given a crossing family $\mathcal{C}$, we say that a set function $f\colon \mathcal{C}\rightarrow\mathbb{R}$ is \emph{crossing-submodular} if 
\begin{equation}
f(U\cap W)+f(U\cup W)\leq f(U)+f(W) \quad \forall U,W\in \mathcal{C}\text{ that cross.} \label{eq:def_crossing_submodular}
\end{equation}
Note that the class of crossing-submodular functions over a fixed family $\mathcal{C}$ is closed under addition.
Classical examples of crossing-submodular functions include, for instance:
\begin{itemize}
    \item \emph{Modular} functions, i.e., functions of the form $U\mapsto \sum_{u\in U} x_u +c$, where $x\in \mathbb{R}^V$ and $c\in\mathbb{R}$.
    \item Functions of the form $U\mapsto c(\delta_D^+(U))$ or $U\mapsto c(\delta_D^-(U))$, defined on a crossing family $\mathcal{C}\subseteq 2^V$, where $D=(V,A)$ is a digraph and $c\in \mathbb{R}_{\geq 0}^A$.
\end{itemize}
In particular, the second example shows that in every digraph $D=(V,A)$, $\{\emptyset\neq U\subsetneq V\colon \delta^-(U)=\emptyset\}$ and $\{\emptyset\neq U\subsetneq V\colon \delta^+(U)=\emptyset\}$ constitute crossing families.

\paragraph{Linear and integer programming.}
Let $A\in\mathbb{Q}^{m\times n}$ and let $b\in \mathbb{Q}^m$.
The system $Ax\leq b$ of linear inequalities is called \emph{totally dual integral (TDI)} if for every cost vector $c\in \mathbb{Z}^n$ for which $\max\{c^\top x\colon Ax\leq b\}$ has an optimum solution, its dual $\min\{ b^\top y\colon A^\top y = c, y\geq 0\}$ has an integral optimum solution. If the right-hand side $b$ is integral and the system $Ax\leq b$ is TDI, then the polyhedron $P\coloneqq \{x\colon Ax\leq b\}$ is integral~\cite{EDMONDS1977185}.

A system $Ax\leq b$ is called \emph{box-TDI} if for every pair of integral vectors $\ell\leq u$, the system $Ax\leq b, \ell\leq x\leq u$ is TDI. In particular, every box-TDI system is also TDI. An example for a box-TDI system is the intersection of two `base systems', as stated formally below.

\begin{theorem}[see \cite{S2003V2}, Theorem 48.9 and \cite{F2011}, \S14.4] \label{thm:box_TDI}
 For $i = 1, 2$, let $\mathcal{C}_i$ be a
crossing family over ground set $V$, let $f_i\colon \mathcal{C}_i\rightarrow\mathbb{Z}$ be a crossing-submodular function, and let $k$ be an integer. Then the system
\[x(V)=k,\, x(U)\leq f_i(U)\quad \forall i\in \{1,2\},\,\forall U\in\mathcal{C}_i\] is box-TDI, and therefore TDI.
\end{theorem}

\section{Proof of \cref{thm:strong_orientation} \label{sec:proof-1}}

Let $G=(V,E)$ be a connected graph, and let $\mathcal{C}$ be a crossing family such that $|\delta_G(U)|\geq 2$ for all $U\in \mathcal{C}$. Let $D=(V,A)$ be an arbitrary orientation of $G$. Our goal is to find a strong re-orientation of $D$ for $\mathcal{C}$. We phrase this as the problem of finding an integral point in the intersection of two submodular flow systems, which also happens to be a generalized set covering system.

Let $\mathcal{C}_1\coloneqq \mathcal{C}$ and $\mathcal{C}_2\coloneqq \{V\setminus U:U\in \mathcal{C}\}$. By definition, $\mathcal{C}_1$ is a crossing family, and it can be readily checked that $\mathcal{C}_2$ is also a crossing family. Moreover, the functions $f_i\colon\mathcal{C}_i\rightarrow\mathbb{R}, U\mapsto |\delta_{D}^+(U)|-1$ for $i\in \{1,2\}$, are crossing-submodular. Consider the system
 \begin{equation}
y(\delta_{D}^+(U))-y(\delta_{D}^-(U))\leq f_i(U)\quad\forall i\in\{1,2\},\,\forall U\in\mathcal{C}_i;\, y\in\mathbb{R}^{A}.\label{eq:intersection_submodular_flows}
 \end{equation} This system can be rewritten as follows:
 \begin{equation}
\sum_{a\in\delta_{D}^-(U)} y_a + \sum_{b\in \delta_{D}^+(U)} (1-y_b)\geq 1\quad\forall i\in\{1,2\},\,\forall U\in\mathcal{C}_i;\, y\in\mathbb{R}^{A}.\label{eq:intersection_GSC}
 \end{equation} We finish the proof in three steps. \begin{enumerate}
 \item $D$ has a strong re-orientation for $\mathcal{C}$ if and only if \eqref{eq:intersection_GSC} admits a $0,1$ solution; let us explain. Given a $0,1$ solution $y^\star$, the digraph obtained from $D$ after flipping the arcs in $\{a\in A:y^\star_a=1\}$ has at least one outgoing arc for each $U\in \mathcal{C}_1\cup \mathcal{C}_2$, implying in turn that the new digraph is strongly connected for $\mathcal{C}$. By reversing the implications we obtain the converse.
 \item If \eqref{eq:intersection_GSC} has an integral solution $\bar{y}$, then it has a $0,1$ solution. To see this, define $y^\star\in \{0,1\}^A$ as follows: for each $a\in A$, $y^\star_a\coloneqq 1$ if $\bar{y}_a\geq 1$, and $y^\star_a\coloneqq 0$ if $\bar{y}_a\leq 0$. Then $y^\star$ is also solution to the system. Suppose for a contradiction that $y^\star$ violates one of the inequalities, say for some $U\in \mathcal{C}_i$. Then $y^\star_a=0$ for all $a\in \delta^-_D(U)$, and $y^\star_b=1$ for all $b\in \delta^+_D(U)$. However, this implies that $\bar{y}_a\leq 0$ for all $a\in \delta^-_D(U)$, and $\bar{y}_b\geq 1$ for all $b\in \delta^+_D(U)$, so $\sum_{a\in\delta_{D}^-(U)} \bar{y}_a + \sum_{b\in \delta_{D}^+(U)} (1-\bar{y}_b)\leq 0$, a contradiction to the feasibility of~$\bar{y}$.
 \item The system \eqref{eq:intersection_submodular_flows} has an integral solution. To see this, let $x'\in\mathbb{R}^V$ be given by $x'_v\coloneqq \frac{1}{2}\cdot(|\delta^+_{D}(v)|-|\delta^-_{D}(v)|)$. Note that $x'(V)=0$. For every $U\in\mathcal{C}_i$, we have
\[
x'(U)
=\frac{1}{2}\cdot (|\delta^+_{D}(U)|-|\delta^-_{D}(U)|) 
= |\delta^+_{D}(U)| - \frac{1}{2}\cdot \underbrace{\left(|\delta^+_{D}(U)|+|\delta^-_{D}(U)|\right)}_{\geq 2}\leq |\delta^+_{D}(U)|-1 
= f_i(U).
\]
Thus, $x'$ is a feasible solution to the system 
\[x(V)=0;\, x(U)\leq f_i(U),\,\forall i\in\{1,2\},\forall U\in\mathcal{C}_i.\]
By \cref{thm:box_TDI}, this system has an integral solution $\bar{x}$. As $D$ is weakly connected and $\bar{x}(V)=0$, there exists an integral $\bar{x}$-transshipment $\bar{y}\in \mathbb{Z}^{A}$, by \cref{b-transshipment}.
As $\bar{y}(\delta_{D}^+(U))-\bar{y}(\delta_{D}^-(U))=\bar{x}(U)$ for all $U\subseteq V$, $\bar{y}$ is an integral solution to \eqref{eq:intersection_submodular_flows}.
 \end{enumerate}
 The last step tells us that \eqref{eq:intersection_submodular_flows}, and therefore \eqref{eq:intersection_GSC}, has an integral solution $\bar{y}$. By the second step, the system also has a $0,1$ solution $y^\star$. By the first step, $D$ has a strong re-orientation for $\mathcal{C}$, as required.
 \qed
 
 \medskip
 
Note that the second step in the proof used crucially that every inequality involved in the system \eqref{eq:intersection_submodular_flows} is a \emph{generalized set covering} inequality, while the third step used crucially that the system is the intersection of two submodular flow systems on a weakly connected digraph; surprisingly, this system is TDI~\cite{Abdi2024}.

\section{Proof of \cref{thm:tree_orientation} \label{sec:proof-2}}

Let $(D=(V,A),w)$ be a weighted digraph with $w\in \{0,1\}^A$, where the weight of every dicut is at least $2$. Suppose the set $A_1$ of the weight-$1$ arcs forms a weakly bridge-connected subdigraph. Our goal is to find two disjoint dijoins of $D$ contained in $A_1$. We proceed by induction on $|V|$, and in three steps. The first two steps are routine, and the third step uses \cref{thm:strong_orientation} to prove the base case.
\begin{enumerate}
\item Suppose $A_1$ contains a cycle $C$. We then contract $C$ in $(D,w)$ to obtain another weighted digraph $(D,w)/C$; the weight of an arc that is not contracted remains the same. Clearly, the minimum weight of a dicut of $(D,w)/C$ remains at least $2$, so by the induction hypothesis, there exist two disjoint dijoins $J_1,J_2$ of $D/C$ contained in the weight-$1$ arcs of $(D,w)/C$. We then uncontract $C$, traverse the cycle in some order, denote by $J'_1\subseteq C$ the set of arcs that agree with this order, and denote $J'_2\coloneqq C\setminus J'_1$. It can be readily checked that $J_1\cup J'_1,J_2\cup J'_2$ are disjoint dijoins of $D$ contained in $A_1$.
\item Suppose $V$ contains a node $v_0$ that is not incident to any weight-$1$ arc. We then delete $v_0$, and add a weight-$0$ arc from every in-neighbor of $v_0$ to every out-neighbor of $v_0$.\footnote{Note that this reduction does not necessarily preserve planarity, and so cannot be applied to drop the `spanning' condition in (\cite{CHUDNOVSKY201618}, Theorem 2.2).} Denote by $(D',w')$ the new weighted digraph. It can be readily checked that the minimum weight of a dicut in $(D',w')$ is at least $2$, so by the induction hypothesis, we find two disjoint dijoins $J_1,J_2$ of $D'$ contained in the weight-$1$ arcs of $(D',w')$. It can be readily checked that $J_1,J_2$ are also disjoint dijoins of $D$ contained in the weight-$1$ arcs of $(D,w)$.
\item We arrive at the base case where $A_1$ contains no cycle, and the weight-$1$ arcs span every vertex of $D$. As $A_1$ forms a weakly bridge-connected subdigraph, it follows that the underlying undirected graph of $D_1\coloneqq(V,A_1)$ is a (spanning) tree. Let $\mathcal{C}:=\{\emptyset\neq U \subsetneq V:\delta^-_D(U)= \emptyset\}$. Since every dicut of $(D,w)$ has weight at least $2$, it follows that $|\delta_{D_1}(U)|\geq 2$ for all $U\in \mathcal{C}$. Thus, by \cref{thm:strong_orientation}, there exists a subset $J\subseteq A_1$ such that every set in $\mathcal{C}$ has at least one outgoing and at least one incoming arc in $(A_1\setminus J)\cup J^{-1}$, where $J^{-1}$ is obtained from $J$ after flipping each arc in the set. Subsequently, $J$ and $A_1\setminus J$ are disjoint dijoins of $D$.\qed
\end{enumerate}

\section{The case of two connected components \label{sec:2_components}}

Given that in the counterexample to the Edmonds-Giles conjecture displayed in \cref{fig:Schrijvers_example}, as well as all other known counterexamples~\cite{Cornuejols02,Williams04}, the weight-$1$ arcs form at least three weakly connected components, it appears natural to ask whether \cref{thm:tree_orientation} extends to the case where the weight-$1$ arcs have at most two weakly bridge-connected components. The same argument again allow us to restrict ourselves to the case where the weight-$1$ arcs form a spanning forest with at most $2$ weakly connected components. 

The crux lies in finding $\bar{y}$.
While we may obtain the vectors $x'$ and $\bar{x}$ via the same reasoning as before, the fact that $D_1$ is no longer weakly connected means that there does not necessarily exist an $\bar{x}$-transshipment in $D_1$ anymore. 
More precisely, such a transshipment exists if and only if we additionally require $\bar{x}(V_1)=\bar{x}(V_2)=0$, where $V_1$ and $V_2$ are the vertex sets of the two weakly connected components of $D_1$. 
Unfortunately, the system
\begin{equation}x(V_1)=x(V_2)=0;\, x(U)\leq f_i(U),\, \forall i\in\{1,2\}, \forall U\in\mathcal{C}_i\label{eq:new_system}\end{equation}
is no longer TDI in general. Even worse, the polytope
\begin{equation}\{x\colon x(V_1)=x(V_2)=0;\, x(U)\leq f_i(U),\, \forall i\in\{1,2\}, \forall U\in\mathcal{C}_i\}\label{eq:polytope_2_comp}\end{equation}
does not even need to be integral; a counterexample is presented in \cref{app:not_integral}.
While this does not rule out the existence of an integral solution to \eqref{eq:new_system}, which would be sufficient for our purposes, obtaining it lies beyond the scope of our techniques.

\section{Strengthening sets for a crossing family}\label{sec:conjecture}

Let $D=(V,A)$ be a digraph, and let $\mathcal{C}$ be a crossing family over ground set $V$. A \emph{strengthening set of $D$ for $\mathcal{C}$} is a subset $J\subseteq A$ such that every set in $\mathcal{C}$ has at least one outgoing arc after the arcs in $J$ are flipped. We make the following conjecture.

\begin{conjecture}\label{strengthening-CN}
Let $\tau\geq 2$ be an integer, let $D=(V,A)$ be a weakly connected digraph, and let $\mathcal{C}$ be a crossing family such that $|\delta_D^-(U)|+(\tau-1)|\delta_D^+(U)|\geq \tau$ for all $U\in \mathcal{C}$. Then $A$ can be partitioned into $\tau$ strengthening sets of $D$ for $\mathcal{C}$.
\end{conjecture}

Observe that the inequality equivalently asks that $|\delta_D(U)|\geq 2$ for all $U\in \mathcal{C}$, and $|\delta^-_D(U)|\geq \tau$ for all $U\in \mathcal{C}$ such that $\delta^+_D(U)=\emptyset$. Both of these inequalities are necessary for $A$ to be partitioned into $\tau$ strengthening sets of $D$ for $\mathcal{C}$.

\cref{strengthening-CN} for $\tau=2$ is equivalent to \Cref{thm:strong_orientation}. It is equivalent to Woodall's conjecture if $\mathcal{C}=\{\emptyset\neq U\subsetneq V : \delta_D^+(U)=\emptyset\}$. 
It becomes (\cite{CHUDNOVSKY201618}, Conjecture 1.5) in the case where $\mathcal{C}=\{\emptyset\neq U\subsetneq V : \delta_{D\cup D'}^+(U)=\emptyset\}$ for another digraph $D'$ over vertex set $V$. Finally, when $\mathcal{C}=\{U\subsetneq V:U\neq \emptyset\}$, we obtain the following conjecture that appears in an unpublished note by Schrijver~\cite{Schrijver-note}.

Given a digraph $D=(V,A)$, a \emph{strengthening set} is a subset $J\subseteq A$ such that $(D\setminus J)\cup J^{-1}$ is strongly connected. 

\begin{conjecture}[\cite{Schrijver-note}]\label{schrijver-CN}
Let $\tau\geq 2$ be an integer, and let $D=(V,A)$ be a digraph where every dicut has size at least $\tau$. Then $A$ can be partitioned into $\tau$ strengthening sets.
\end{conjecture}

It can be readily checked that if $A$ is partitioned into $\tau$ strengthening sets, then it has $\tau$ disjoint dijoins. Conversely, there is a digraph $D'$ such that every packing of $\tau$ dijoins of $D'$ yields a partition of $A$ into $\tau$ strengthening sets of $D$~\cite{Schrijver-note}.\footnote{To obtain $D'$, simply replace every arc $a=(u,v)$ of $D$ with $(u,t_a)$ and $(\tau-1)$ parallel copies of $(v,t_a)$ for a new node $t_a$.} In summary, for a specific digraph, the conjecture above is at least as strong as Woodall's conjecture, while for all digraphs, the two conjectures are equivalent.

\section*{Acknowledgements}

We would like to thank G\'{e}rard Cornu\'{e}jols and Giacomo Zambelli for fruitful discussions about this work. This work was supported in part by EPSRC grant EP/X030989/1.

\paragraph{Data Availability Statement.} No data are associated with this article. Data sharing is not applicable to this article.

\bibliographystyle{plainurl}
\bibliography{Woodalls_conjecture}

\begin{thebibliography}{10}

\bibitem{Abdi2024}
Ahmad Abdi, G{\'e}rard Cornu{\'e}jols, and Giacomo Zambelli.
\newblock Arc connectivity and submodular flows in digraphs.
\newblock {\em Combinatorica}, 2024.
\newblock \href {https://doi.org/10.1007/s00493-024-00108-0}
  {\path{doi:10.1007/s00493-024-00108-0}}.

\bibitem{Abdi2023}
Ahmad Abdi, G\'{e}rard Cornu\'{e}jols, and Michael Zlatin.
\newblock On packing dijoins in digraphs and weighted digraphs.
\newblock {\em SIAM Journal on Discrete Mathematics}, 37(4):2417--2461, 2023.
\newblock \href {https://doi.org/10.1137/22M1506511}
  {\path{doi:10.1137/22M1506511}}.

\bibitem{CHUDNOVSKY201618}
Maria Chudnovsky, Katherine Edwards, Ringi Kim, Alex Scott, and Paul Seymour.
\newblock Disjoint dijoins.
\newblock {\em Journal of Combinatorial Theory, Series B}, 120:18--35, 2016.
\newblock \href {https://doi.org/10.1016/j.jctb.2016.04.002}
  {\path{doi:10.1016/j.jctb.2016.04.002}}.

\bibitem{Cornuejols02}
G{\'e}rard Cornu{\'e}jols and Bertrand Guenin.
\newblock On dijoins.
\newblock {\em Discrete Mathematics}, 243(1):213--216, 2002.
\newblock \href {https://doi.org/10.1016/S0012-365X(01)00209-6}
  {\path{doi:10.1016/S0012-365X(01)00209-6}}.

\bibitem{10.1007/978-3-031-59835-7_6}
G{\'e}rard Cornu{\'e}jols, Siyue Liu, and R.~Ravi.
\newblock Approximately packing dijoins via nowhere-zero flows.
\newblock In Jens Vygen and Jaros{\l}aw Byrka, editors, {\em Integer
  Programming and Combinatorial Optimization}. Springer, 2024.
\newblock \href {https://doi.org/10.1007/978-3-031-59835-7_6}
  {\path{doi:10.1007/978-3-031-59835-7_6}}.

\bibitem{EDMONDS1977185}
Jack Edmonds and Rick Giles.
\newblock A min-max relation for submodular functions on graphs.
\newblock In P.L. Hammer, E.L. Johnson, B.H. Korte, and G.L. Nemhauser,
  editors, {\em Studies in Integer Programming}, volume~1 of {\em Annals of
  Discrete Mathematics}, pages 185--204. Elsevier, 1977.
\newblock \href {https://doi.org/10.1016/S0167-5060(08)70734-9}
  {\path{doi:10.1016/S0167-5060(08)70734-9}}.

\bibitem{F2011}
Andr{\'a}s Frank.
\newblock {\em Connections in Combinatorial Optimization}.
\newblock Oxford Lecture Series in Mathematics and Its Applications. Oxford
  University Press, 2011.
\newblock \href {https://doi.org/10.1016/j.dam.2011.09.003}
  {\path{doi:10.1016/j.dam.2011.09.003}}.

\bibitem{Schrijver-note}
Alexander Schrijver.
\newblock Observations on {W}oodall's conjecture.
\newblock Accessed online at
  \url{https://homepages.cwi.nl/~lex/files/woodall.pdf}.

\bibitem{SCHRIJVER1980213}
Alexander Schrijver.
\newblock A counterexample to a conjecture of {E}dmonds and {G}iles.
\newblock {\em Discrete Mathematics}, 32(2):213--214, 1980.
\newblock \href {https://doi.org/10.1016/0012-365X(80)90057-6}
  {\path{doi:10.1016/0012-365X(80)90057-6}}.

\bibitem{S2003V1}
Alexander Schrijver.
\newblock {\em Combinatorial Optimization: Polyhedra and Efficiency Volume 1}.
\newblock Springer Science \& Business Media, 2003.

\bibitem{S2003V2}
Alexander Schrijver.
\newblock {\em Combinatorial Optimization: Polyhedra and Efficiency Volume 2}.
\newblock Springer Science \& Business Media, 2003.

\bibitem{Shepherd05}
Bruce Shepherd and Adrian Vetta.
\newblock Visualizing, finding and packing dijoins.
\newblock In D.~Avis, A.~Hertz, and O.~Marcotte, editors, {\em Graph Theory and
  Combinatorial Optimization}, chapter~8, pages 219--254. Springer Verlag.
\newblock \href {https://doi.org/10.1007/0-387-25592-3_8}
  {\path{doi:10.1007/0-387-25592-3_8}}.

\bibitem{Williams04}
Aaron Williams.
\newblock Packing directed joins.
\newblock Master's thesis, University of Waterloo, 2004.
\newblock URL: \url{http://hdl.handle.net/10012/1024}.

\bibitem{Woodall}
Douglas~R. Woodall.
\newblock Menger and {K}{\"o}nig systems.
\newblock In Yousef Alavi and Don~R. Lick, editors, {\em Theory and
  Applications of Graphs}, pages 620--635. Springer, 1978.
\newblock \href {https://doi.org/10.1007/BFb0070416}
  {\path{doi:10.1007/BFb0070416}}.

\end{thebibliography}

\FloatBarrier

\appendix

\section{An example where \eqref{eq:polytope_2_comp} is not integral.\label{app:not_integral}}

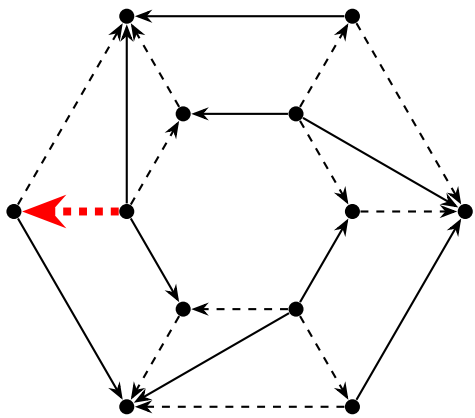
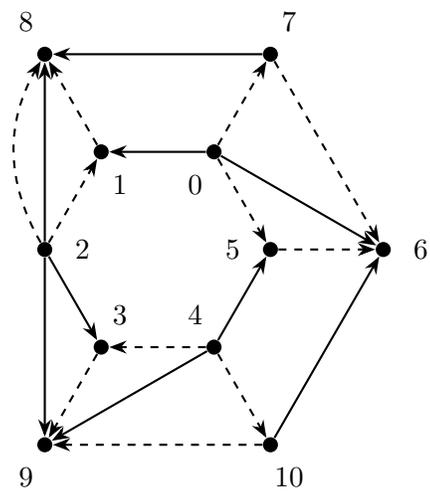
\begin{figure}
\begin{subfigure}[t]{0.45\textwidth}
\centering
\begin{tikzpicture}
\foreach \i in {1,...,6}
{
  \node[circle, fill, inner sep = 0pt, minimum size = 2mm] (A\i) at ({1.5*cos(\i*60)},{1.5*sin(\i*60)}) {};
  \node[circle, fill, inner sep = 0pt, minimum size = 2mm] (B\i) at ({3*cos(\i*60)},{3*sin(\i*60)}) {};
}
\foreach \i in {1,2,4,5,6}
{
    \draw[thick, dashed, -{Stealth}] (A\i)--(B\i);
}
\draw[line width=3pt, red, dashed, -{Stealth}] (A3)--(B3);
\foreach \i/\j in {1/6,3/2,5/4}
{
    \draw[thick, dashed,-{Stealth}] (A\i)--(A\j);
    \draw[thick,-{Stealth}] (A\i)--(B\j);
    \draw[thick, dashed, -{Stealth}] (B\i)--(B\j);
}
\foreach \i/\j in {5/6,3/4,1/2}
{
    \draw[thick,-{Stealth}] (A\i)--(A\j);
    \draw[thick,-{Stealth}] (B\i)--(B\j);
}
\foreach \i in {6,...,11}
{
\node at ({3.5*cos((\i)*60)},{3.5*sin((\i)*60)}) {\phantom{\i}};
}
\end{tikzpicture}
\subcaption{Schrijver's example. The thick dashed (red) arc of weight $0$ is contracted to obtain the example on the right.\label{fig:example_non_integral_two_comp_a}}
\end{subfigure}
\hspace{1cm}
\begin{subfigure}[t]{0.45\textwidth}
\centering
\begin{tikzpicture}
\foreach \i in {1,2,4,5,6}
{
  \node[circle, fill, inner sep = 0pt, minimum size = 2mm] (A\i) at ({1.5*cos(\i*60)},{1.5*sin(\i*60)}) {};
  \node[circle, fill, inner sep = 0pt, minimum size = 2mm] (B\i) at ({3*cos(\i*60)},{3*sin(\i*60)}) {};
  \draw[thick, dashed, -{Stealth}] (A\i)--(B\i);
}
\node[circle, fill, inner sep = 0pt, minimum size = 2mm] (A3) at ({1.5*cos(3*60)},{1.5*sin(3*60)}) {};
\foreach \i/\j in {1/6,5/4}
{
    \draw[thick, dashed,-{Stealth}] (A\i)--(A\j);
    \draw[thick,-{Stealth}] (A\i)--(B\j);
    \draw[thick, dashed, -{Stealth}] (B\i)--(B\j);
}
\draw[thick, dashed,-{Stealth}] (A3)--(A2);
\draw[thick,-{Stealth}] (A3)--(B2);
\draw[thick, dashed, -{Stealth}] (A3)to[bend left = 30] (B2);
\foreach \i/\j in {5/6,1/2}
{
    \draw[thick,-{Stealth}] (A\i)--(A\j);
    \draw[thick,-{Stealth}] (B\i)--(B\j);
}
\draw[thick,-{Stealth}] (A3)--(A4);
\draw[thick,-{Stealth}] (A3)--(B4);
\foreach \i in {0,...,5}
{
\node at ({cos((\i+1)*60)},{sin((\i+1)*60)}) {\i};
}
\foreach \i in {6,7,8}
{
\node at ({3.5*cos((\i)*60)},{3.5*sin((\i)*60)}) {\i};
}
\foreach \i in {9,10}
{
\node at ({3.5*cos((\i+1)*60)},{3.5*sin((\i+1)*60)}) {\i};
}
\end{tikzpicture}
\subcaption{A weighted digraph where the weight-$1$ arcs form two connected components and every dicut has weight at least $2$.\label{fig:example_non_integral_two_comp_b}}
\end{subfigure}
 \caption{Schrijver's example (left) and an example where \eqref{eq:polytope_2_comp} is not integral (right). Weight-$0$/$1$ arcs are indicated by dashed/solid lines. For convenience, the vertices are numbered in the right example.}   \label{fig:example_non_integral_two_comp}
\end{figure}

\cref{fig:example_non_integral_two_comp_b} shows a weighted digraph $(D,w)$ for which \eqref{eq:polytope_2_comp} has fractional vertices. It is obtained from Schriver's example by contracting one of the weight-$0$ arcs (see \cref{fig:example_non_integral_two_comp_a}).
One example of a fractional vertex is given by 
\[x^\star\coloneqq [ 1, -0.5,  1.5,  0 ,  1 , -1 , -1 ,  0.5, -1, -1,  0.5],\]
where the first entry corresponds to vertex $0$, the second one corresponds to vertex $1$, etc. Feasibility of $x^\star$ can be checked, e.g., via a computer-assisted brute force enumeration of all subsets of the vertex set.

To certify that $x^\star$ constitutes a vertex of \eqref{eq:polytope_2_comp}, we provide the following list of $11$ linearly independent tight constraints:
\begin{align*}
x^\star(\{0\}) &= f_1(\{0\})=2-1 =1\\
x^\star(\{4\}) &= f_1(\{4\})=2-1 =1\\
x^\star(\{6\}) &= f_2(\{6\}) = -1\\
x^\star(\{8\}) &= f_2(\{8\}) = -1\\
x^\star(\{9\}) &= f_2(\{9\}) = -1\\
x^\star(\{3,9\}) &= f_2(\{3,9\}) = -1\\
x^\star(\{1,2,3,8,9\})& = f_2(\{1,2,3,8,9\}) = -1\\
x^\star(\{2,3,4,9,10\}&= f_1(\{2,3,4,9,10\} = 3-1 = 2\\
x^\star(\{0,1,2,3,4,7,8,9,10\} &= f_1(\{0,1,2,3,4,7,8,9,10\}) = 3-1 = 2\\
x^\star(\{0,1,6,10\})&=x^\star(V_1)=0\\
x^\star(\{2,3,4,5,7,8,9\})&=x^\star(V_2)=0\\
\end{align*}
\end{document}